\documentclass[a4paper]{article}
\usepackage{amsmath, amssymb, eucal, amscd}

\newtheorem{thm}{Theorem}[section]
\newtheorem{lem}[thm]{Lemma}
\newtheorem{eg}[thm]{Example}
\newtheorem{prop}[thm]{Proposition}
\newtheorem{cor}[thm]{Corollary}
\newtheorem{defn}[thm]{Definition}
\newtheorem{rem}[thm]{Remark}
\newtheorem{ntn}[thm]{Notation}
\newtheorem{rem-eg}[thm]{Remark and Example}

\newcommand{\smnoind}{{\smallskip\noindent}}

\newcommand{\id}{{\rm id}}

\newcommand{\ti}{\tilde}

\newcommand{\cb}{{\rm cb}}
\newcommand{\CB}{{\rm CB}}
\newcommand{\abs}[1]{\left\vert#1\right\vert}

\newcommand{\sph}{\mathfrak{S}_1}
\newcommand{\CS}{\mathcal{S}}

\newcommand{\CL}{\mathcal{L}}
\newcommand{\OV}{\mathbf{V}} 
\newcommand{\OW}{\mathbf{W}}
\newcommand{\OZ}{\mathbf{Z}}
\newcommand{\BK}{\mathbf{K}}

\newenvironment{prf}{{\noindent \textbf{Proof:} }}{\hfill $\Box$\medskip}

\begin{document}

\title{An abstract characterization of unital operator spaces\footnote{This work is jointly supported by
Hong Kong RGC Research Grant (2160255), the National Natural Science Foundation of China (10771106) and NCET-05-0219}}
\author{Xu-Jian Huang and Chi-Keung Ng}

\maketitle
\begin{abstract}
In this article, we give an abstract characterization of the ``identity'' of an operator space $\OV$ by looking at a quantity $n_\cb(\OV,u)$ which is defined in analogue to a well-known quantity in Banach space theory. 
More precisely, we show that there exists a complete isometry from $\OV$ to some $\CL(H)$ sending $u$ to $\id_H$ if and only if $n_\cb(\OV, u) =1$. 
We will use it to give an abstract characterization of operator systems. 
Moreover, we will show that if $\OV$ is a unital operator space and $\OW$ is a proper complete $M$-ideal, then $\OV/\OW$ is also a unital operator space. 
As a consequece, the quotient of an operator system by a proper complete $M$-ideal is again an operator system. 
In the appendix, we will also give an abstract characterisation of ``non-unital operator systems'' using an idea arose from the definition of $n_\cb(\OV,u)$. 
\end{abstract}

\section{Introduction}

\bigskip

Operator spaces are subspaces of $\CL(H)$ (where $H$ is a Hilbert spaces) together with the induced ``matrix norm structures''. 
The theory of operator spaces is a very important tool in the study of operator algebras. 
Since the discovery of an abstract characterization of operator spaces by Ruan (see e.g. \cite{Ruan-char} or \cite{ER}), there have been many more applications of operator spaces to other branches in functional analysis (see e.g. \cite{ARS}, \cite{IS}, \cite{KR}, \cite{LNR}, \cite{Ng-cohom}, \cite{Ng-ext}, \cite{Ruan-op-amen}, \cite{Ruan-amen-Hopf} \cite{Runde-os-aha} and \cite{Spr}). 

\medskip

In the theory of operator spaces, sometimes the starting point is not general subspaces of $\CL(H)$ but unital subspaces of $\CL(H)$ (i.e. subspaces that contain the identity; see e.g. \cite{BL}). 
It is natural to ask whether there is an abstract characterization of unital operator spaces. 

\medskip

In fact, there is a Banach space characterization for the identity of a unital $C^*$-algebra (see e.g. \cite[Theorem 2]{AW}, \cite[4.1]{BK} and \cite[9.5.16]{Pal}) which gives rise to the concept of ``geometric unitary'' (see e.g. \cite{LNW}). 
More precisely, a norm one element $u$ in a Banach space $X$ is a geometric unitary if certain quantity $n(X;u)$ is non-zero. 
This quantity is also related to numerical indexes of Banach spaces (see e.g. \cite{BD}). 
It is natural to think that certain operator space analogue of such quantity may give a characterization of the ``identity'' of an operator space. 

\medskip

In this paper, we will define and study such a quantity $n_\cb(\OV;u)$ for an operator space $\OV$ and a norm one element $u\in \OV$. 
We will give some properties about it and show that there exists a complete isometry (respective, complete isomorphism) from $\OV$ to some $\CL(H)$ that sends $u$ to $\id_H$ if and only if $n_\cb(\OV;u) = 1$ (respectively, $n_\cb(\OV;u) > 0$). 
We will also give some relations between $n_\cb(\OV;u)$ and $n(V;u)$. 
As an application, we give an abstract characterization for operator systems. 
Furthermore, we will show that if $\OV$ is a unital operator space (respectively, operator system), then any proper complete $M$-ideal $\OW$ of $\OV$ will not contain the identity and the quotient $\OV/\OW$ is again a unital operator space (respectively, operator system). 

\medskip

In the appendix, we give an abstract characterisation for ``non-unital operator systems'' using the similar idea as that for ordinary operator systems as in the main body of the article. 

\medskip

In the following, let us first recall the notion of ``geometric unitary'' for normed spaces. 
Suppose that $X$ is a normed space and $\sph(X)$ is the unit sphere of $X$. 
For any $u\in \sph(X)$, we set
$$\CS(X; u)\ :=\ \{f\in X^*: \|f\| = 1 = f(u)\},$$ 
$$\gamma^{u}(x)\ :=\ \sup\left\{ \|f(x)\|: f\in \CS(X; u)\right\} \qquad (x\in X),$$
as well as
$$n(X;u)\ :=\ \inf\left\{\gamma^{u}(x): x\in \sph(V)\right\}.$$

\medskip

A norm one element $u$ in a normed space $X$ is called a \emph{geometric unitary} (respectively, \emph{strict geometric unitary}) if $n(V;u) > 0$ (respectively, $n(V;u) =1$). 

\medskip

The following result is probably well-known. 
Part (a) of which is one of the motivations behind this work. 
Since the arguments for it is straight forward (and also follows from a similar arguments as those for Theorem \ref{ncb=1} and Proposition \ref{map-ncb}(b)), we leave it to the readers to check them. 

\medskip

\begin{prop}
\label{map-n}
Let $X$ be a normed space and $u\in \sph(X)$. 

\smnoind
(a) $u$ is a strict geometric unitary (respectively, geometric unitary) if and only if there exist a compact Hausdorff space $\Omega$ and an isometry (respectively, a contractive topological injection) $\varphi: X \rightarrow C(\Omega)$ such that $\varphi (u) = 1_\Omega$. 

\smnoind
(b) Suppose that $Y$ is another normed space and $\Psi:X \rightarrow Y$ is a contractive topological injection. 
If $\Psi(u)\in \sph(Y)$, then $n(Y; \Psi(u)) \leq \|\Psi^{-1}\|\cdot n(X; u)$. 
\end{prop}

\bigskip

\section{Complete geometric unitaries}

\bigskip

In analogue to the notion of geometric unitaries, one can define complete geometric unitaries for operator spaces. 
In fact, we will go a step further and start with matricially normed spaces (as defined in \cite[p.217]{Ruan-char}) instead of operator spaces. 

\medskip

\begin{ntn}
Throughout this article, unless specified, $\OV$ is a matricially normed space (i.e. a matrix normed space that satisfies only condition (M2) in \cite[p.20]{ER}) while $V$ is the underlying normed space (or underlying vector space) of $\OV$. 
We will denote by $\sph(\OV)$ the unit sphere of $\OV$. 
\end{ntn}

\medskip

\begin{rem}
\label{rem-mmns}
If $M_\infty(V)$ denote the set of all infinite matrices on $V$ with only finitely many non-zero entries, then by \cite[2.1]{Ruan-char}, one can turn $M_\infty(V)$ into an essential normed $M_\infty(\mathbb{C})$-bimodule (where $M_\infty(\mathbb{C})$ is regarded as a normed subalgebra of $\CL(\ell_2)$). 
Therefore, one can construct the so-called \emph{regularization} of $M_\infty(V)$ (as in \cite[1.2]{Ng-reg}) which produces an operator space $\OV_0$ (by \cite[2.1]{Ng-reg}). 
In fact, $\OV_0$ is the image of $\OV$ in $\OV^{**}$ ($\OV^{**}$ is as defined in \cite[p.41]{ER} rather than \cite{Ruan-char}) together with the induced norm. 
Now, by a similar argument as that of \cite[2.5]{Ng-reg}, if $\OW$ is an operator space, then the canonical map from $\CB(\OV_0, \OW)$ to $\CB(\OV, \OW)$ given by the canonical map $\kappa_V: \OV \rightarrow \OV_0$ is a complete isometry. 
\end{rem}

\medskip

\subsection{Definition and main results}

\medskip

Let $u\in \sph(\OV)$. 
For any $n\in \mathbb{N}$, we define 
$$\CS_n(\OV;u)\ :=\ \{\varphi\in \CB(\OV, M_n): \|\varphi\|_\cb \leq 1; \varphi (u) = I_n\},$$ 
$$\gamma_k^{u}(x)\ :=\ \sup\left\{ \|\varphi_k(x)\|: \varphi\in \CS_n(\OV;u); n\in\mathbb{N}\right\} \qquad (k\in \mathbb{N}; x\in M_k (\OV))$$
as well as
$$n_\cb(\OV; u)\ :=\ \inf\left\{\gamma_k^{u}(x): x\in \sph(M_k(\OV)); k\in \mathbb{N}\right\}.$$

\medskip

In the case when $\OV$ is an operator algebra and $u$ being the identity of $\OV$, the quantity $\gamma_k^{u}(x)$ was defined in \cite[p.192]{BRS} and is called the \emph{$k$th numerical radius} of $x\in \OV$. 

\medskip

\begin{defn}
A norm one element $u$ in $\OV$ is called a \emph{complete geometric unitary} (respectively, \emph{complete strict geometric unitary}) if $n_\cb(\OV; u) > 0$ (respectively, $n_\cb(\OV;u) =1$). 
\end{defn}

\medskip

\begin{rem}
\label{rem-Def}
(a) $\CS_1(\OV;u) = \CS(V;u)$. 

\smnoind
(b) If $\OV_0$ is the operator space as given in Remark \ref{rem-mmns}, then $\CS_n(\OV, u) = \CS_n(\OV_0, \kappa_V(u))$. 

\smnoind
(c) $\CS_n(\OV;u)$ is a compact convex set under the point-norm topology. 
In fact, it is not hard to see that $\CS_n(\OV;u)$ is point-norm closed in the closed unit ball of $\CB(\OV,M_n)$ which is compact under the point-norm topology (by part (b) and the corresponding fact for operator spaces). 

\smnoind
(d) Suppose that $\OV$ is a matrix normed subspace of a matricially normed space $\OW$ and $\gamma_k^{u,W}$ is defined by $\CS_n(\OW, u)$ in a similar way as $\gamma_k^u$. 
By the Arveson extension theorem and part (b), the canonical map gives a surjection from $\CS(\OW; u)$ to $\CS(\OV; u)$. 
Therefore, $\gamma_k^u = \gamma_k^{u,W}$ for all $k\in \mathbb{N}$ (so we can use $\gamma_k^u$ instead of $\gamma_k^{u,V}$). 
Consequently, $n_\cb(\OW,u) \leq n_\cb(\OV,u)$. 

\smnoind
(e) Suppose that $X$ is a vector space and $\sigma_k$ is a seminorm on $M_k(X)$ ($k\in \mathbb{N}$) satisfying the following two conditions: 
for any $m,n\in \mathbb{N}$, $x\in M_m(X)$, $y\in M_n(Y)$, $\alpha\in M_{n,m}$ and $\beta\in M_{m,n}$,
\begin{enumerate}
\item [\rm I.] $\sigma_{m+n}(x\oplus y) \leq \max \{ \sigma_m(x), \sigma_n(y)\}$; 
\item [\rm II.] $\sigma_n(\alpha\cdot x \cdot \beta) \leq \|\alpha\|\sigma_m(x)\|\beta\|$.
\end{enumerate}
Let $N = \{v\in V: \sigma_1^u(v) = 0\}$. 
Then $\sigma_k$ induces a semi-norm $\ti \sigma_k$ on $M_k(X/N)$ (since $\sigma_k((x_{i,j})_{i,j=1}^k) \leq \sum_{i,j=1}^n \sigma_1(x_{i,j})$) which also satisfies the two conditions as the above. 
Therefore, by \cite[2.3.6]{ER}, $\{\ti\sigma_k\}$ gives an operator space structure on $X/N$. 
\end{rem}

\medskip

\begin{lem}
\label{def-Vu}
$\gamma_k^u$ is a semi-norm on $M_k(\OV)$ ($k\in\mathbb{N}$) satisfying the two properties in Remark \ref{rem-Def}(e) and they induce an operator space structure on $\OV_u := \OV / N_u$ where $N_u := \{v\in \OV: \gamma_1^u(v) = 0\}$. 
\end{lem}

\medskip

In the following, when we talk about the operator space $\OV_u$, we consider the operator space structure as given in the above lemma (even when $N_u = (0)$). 
Moreover, we denote by $Q_u$ the canonical complete contraction from $\OV$ to $\OV_u$. 

\medskip

\begin{lem}
\label{equiv-ncb>0}
Let $\OV$ be a matricially normed space and $u\in \sph(\OV)$. 

\smnoind
(a) $n_\cb(\OV; u) > 0$ if and only if $Q_u$ is a complete isomorphism. 
In this case, $n_\cb(\OV; u) = \|Q_u^{-1}\|_\cb^{-1}$. 

\smnoind
(b) $n_{cb}(\OV_u, Q_u(u)) = 1$. 
\end{lem}
\begin{prf}
(a) Since the norm on $M_n(\OV_u)$ is given by $\ti\gamma_n^u((Q_u)_n(x)) := \gamma_n^u(x)$ ($n\in\mathbb{N}, x\in M_n(\OV)$) and 
\begin{equation}
\label{equ-ncb}
n_\cb(\OV; u) 
\ = \ \sup\ \{ \lambda \in \mathbb{R}_+: \lambda \|x\| \leq \gamma_n^u(x); n\in \mathbb{N}; z\in M_n(\OV)\},
\end{equation}
we obtain the first statement. 
If $n_\cb(\OV;u) > 0$, then 
$$n_\cb(\OV; u)^{-1} 
\ = \ \inf\ \{ \mu \in \mathbb{R}_+: \|(Q_u^{-1})_n(y)\| \leq \mu \ti \gamma_n^u(y); n\in \mathbb{N}; y\in M_n(\OV_u)\}$$
which gives the second statement.

\smnoind
(b) Consider the map $(\ti Q_u)_n: \CB(\OV_u, M_n) \rightarrow \CB(\OV,M_n)$ given by composition with $Q_u$. 
It is clear that $(\ti Q_u)_n(\CS_n(\OV_u; Q_u(u))) \subseteq \CS_n(\OV;u)$ ($n\in \mathbb{N}$). 
On the other hand, for any $\varphi \in \CS_n(\OV;u)$, we have 
$$\|\varphi_k(x)\| \leq \gamma^u_k(x) = \ti \gamma^u_k((Q_u)_k(x)) \qquad (x\in M_k(\OV)).$$
Hence there exists $\psi\in \CB(\OV_u, M_n)$ with $\varphi = \psi \circ Q_u$ and $\|\psi\|_\cb \leq 1$. 
This shows that $(\ti Q_u)_n$ is a surjection (and hence a bijection) from $\CS_n(\OV_u; Q_u(u))$ to $\CS_n(\OV; u)$. 
Consequently, $N_{Q_u(u)} = (0)$ and $\ti \gamma^u_k = \gamma^{Q_u(u)}_k$ on $M_k(V_u) = M_k(V_{Q_u(u)})$. 
Hence $Q_{Q_u(u)}$ is a completely isometric isomorphism and $n_\cb(\OV_u; Q_u(u)) = 1$ by part (a). 
\end{prf}

\medskip

\begin{prop}
\label{map-ncb}
Let $\OV$ and $\OW$ be matricially normed spaces.
Suppose that $\Psi:\OV \rightarrow \OW$ is a complete contraction and $u\in \sph(\OV)$ such that $\Psi(u)\in \sph(\OW)$.

\smnoind
(a) There exists a complete contraction $\Psi_u: \OV_u \rightarrow \OW_{\Psi(u)}$ with $Q_{\Psi(u)}\circ \Psi = \Psi_u \circ Q_u$. 

\smnoind
(b) If $\Psi$ is a complete topological injection, then $n_\cb(\OW; \Psi(u)) \leq \|\Psi^{-1}\|_\cb\cdot n_\cb(\OV; u)$. 
\end{prop}
\begin{prf}
(a) The composition with $\Psi$ gives a map $\ti \Psi_n: \CS_n(\OW;\Psi(u)) \rightarrow \CS_n(\OV;u)$ and we have $\gamma^{\Psi(u)}_n\circ \Psi \leq \gamma^u_n$ ($n\in \mathbb{N}$). 
This induces the required map $\Psi_u$. 

\smnoind
(b) By Remark \ref{rem-Def}(d), Equality \eqref{equ-ncb} and the inequality in the argument for part (a), 
\begin{eqnarray*}
n_\cb(\OW;\Psi(u))
& \leq & n_\cb(\Psi(\OV);\Psi(u))\\
& = & \sup \{ \lambda\in \mathbb{R}_+: \lambda \|\Psi_n(x)\| \leq \gamma^{\Psi(u)}_n(\Psi_n(x)); n\in \mathbb{N}; x\in M_n(\OV)\}\\ 
& \leq & \sup \{ \lambda\in \mathbb{R}_+: \lambda \|\Psi^{-1}\|_\cb^{-1}\|x\| \leq \gamma^{u}_n(x); n\in \mathbb{N}; x\in M_n(\OV)\}\\
& = & \|\Psi^{-1}\|_\cb\cdot n_\cb(\OV;u).
\end{eqnarray*}
\end{prf} 

\medskip

Consequently, for any complete contraction $\Psi: \OV\rightarrow \CL(H)$ with $\Psi(u) = \id_H$ (where $u\in \sph(\OV)$), there exists a complete contraction $\Psi_u: \OV_u \rightarrow \CL(H)$ with $\Psi = \Psi_u \circ Q_u$. 

\medskip

The following is our first theorem which tells us that one can use strict complete geometric unitary to describe the ``identity'' of an operator space. 
Note that one direction of this theorem appeared in the disguised form in \cite[Proposition 1.5]{BRS} in the case of operator algebras and its proof can be carried over directly to the case of operator spaces. 
However, we will give its easy proof here for completeness. 
To do this, we need the following well-known fact (an argument for it can be found in the proof of \cite[Proposition 1.5]{BRS}). 

\medskip 

\begin{lem}
\label{supPTP}
For any Hilbert space $H$ and any $T \in \CL(H^n)$, 
$$\|T\| \ = \ \sup \left\{ \left\|(P\otimes 1)T(P\otimes 1)\right\|: P\in \CL(H) {\rm\ is\ a\ finite\ dimensional\ projection}\right\}.$$
\end{lem}

\medskip

\begin{thm}
\label{ncb=1}
Let $\OV$ be a matricially normed space and $u\in \sph(\OV)$. 
Then $u$ is a complete strict geometric unitary of $\OV$ if and only if there exists a Hilbert space $H$ and a complete isometry $\Theta: \OV\rightarrow \CL(H)$ such that $\Theta(u) = \id_H$. 
\end{thm}
\begin{prf}
$\Leftarrow$). 
For any finite dimensional projection $P$ on $H$, if the rank of $P$ is $n$, then $x \mapsto P\Theta(x)P$ can be regarded as an element of $\CS_n(\OV;u)$. 
Therefore, $\|\cdot\|_k \leq \gamma^u_k$ (and hence equal) on $M_k(V)$ ($k\in \mathbb{N}$) by Lemma \ref{supPTP}. 
Thus, $Q_u$ is a complete isometry and $n_\cb(\OV, u) = 1$ (see e.g. Lemma \ref{equiv-ncb>0}(a)). 

\smnoind
$\Rightarrow$). 
If $\Omega_n$ is the compact Hausdorff space $\CS_n(\OV;u)$ (under the point-norm topology; see Remark \ref{rem-Def}(c)), then $A = \bigoplus_{n=1}^\infty C(\Omega_n, M_n)$ is a unital $C^*$-algebra. 
We define $\Theta: \OV \rightarrow A$ by $\Theta(x) = (\Theta^{(n)}(x))_{n\in \mathbb{N}}$ where $\Theta^{(n)}: \OV \rightarrow C(\Omega_n, M_n)$ given by
$$\Theta^{(n)}(x)(\varphi)\ =\ \varphi(x) \qquad (x\in \OV; \varphi\in \Omega_n).$$ 
For any $k\in \mathbb{N}$ and $z\in M_k(V)$, we have $\Theta_k(z) \in M_k(A) \cong \bigoplus_{n=1}^\infty C(\Omega_n, M_{nk})$ and 
$\Theta^{(n)}_k(z)(\varphi) = \varphi_k(z)$ ($\varphi\in \Omega_n$). 
Thus, 
$$\|\Theta_k(z)\|\ =\ \sup \{ \|\varphi_k(z)\|: n\in \mathbb{N}; \varphi\in \CS_n(\OV;u)\}\ =\ \gamma^u_k(z) \ = \ \|z\|$$
(as $n_\cb(\OV, u) =1$) and so, $\Theta$ is a complete isometry. 
It is not hard to see that $\Theta(u) = \id_H$. 
\end{prf}

\medskip

One can also prove the necessity of the above theorem by taking $H = \bigoplus_{n=1}^\infty \bigoplus_{\varphi\in \CS_n(\OV;u)} \mathbb{C}^n$ and adapting the argument of \cite[2.3.5]{ER}. 

\medskip

\begin{cor}
\label{ncb+>ci}
(a) Suppose that $\OV$ is a matricially normed space. 
Then $n_\cb(\OV;u) > 0$ if and only if there exists a Hilbert space $H$ and a completely contractive complete topological injection $\Psi: \OV\rightarrow \CL(H)$ such that $\Psi(u) = \id_H$. 

\smnoind
(b) If $\OV$ is an operator space, then $n_\cb(\OV;u) = n_\cb(\OV^{**};u)$. 
\end{cor}
\begin{prf}
(a) If $n_\cb(\OV;u) > 0$, then by Lemma \ref{equiv-ncb>0}(a), $Q_u: \OV \rightarrow \OV_u$ is a completely contractive complete isomorphism, and one can apply Theorem \ref{ncb=1} to $\OV_u$ (because of Lemma \ref{equiv-ncb>0}(b)). 
The converse follows from Proposition \ref{map-ncb}(b) and Theorem \ref{ncb=1}. 

\smnoind
(b) Note that if $n_\cb(\OV,u) = 0$, then $n_\cb(\OV^{**}, u) = 0$ (by Remark \ref{rem-Def}(d)). 
Moreover, if $n_\cb(\OV, u) = 1$, then there is a Hilbert space $H$ and a complete isometry $\Theta: \OV \rightarrow \CL(H)$ such that $\Theta(u) =\id_H$ (Theorem \ref{ncb=1}), and so, $\Theta^{**}: \OV^{**}\rightarrow \CL(H)^{**}$ is a complete isometry with $\Theta^{**}(u) = \id_H$ which implies that $n_\cb(\OV^{**}, u) =1$ (note that there exists a faithful unital representation for $\CL(H)^{**}$ and we can apply Theorem \ref{ncb=1}). 
We now suppose that $n_\cb(\OV,u) > 0$. 
Then by Lemma \ref{equiv-ncb>0}(a), $Q_u$ is a complete isomorphism and so is $Q_u^{**}: \OV^{**} \rightarrow (\OV_u)^{**}$. 
By Lemma \ref{equiv-ncb>0}(b), we have $n_\cb(\OV_u, Q_u(u)) = 1$ and the above implies that $n_\cb((\OV_u)^{**}, Q_u(u)) = 1$. 
Consequently, Lemma \ref{equiv-ncb>0}(a), Proposition \ref{map-ncb}(b) as well as Remark \ref{rem-Def}(d) tell us that 
$$n_\cb(\OV;u)\ = \ \|Q_u^{-1}\|_\cb^{-1} \ = \ \|(Q_u^{**})^{-1}\|_\cb^{-1}\ \leq \ n_\cb(\OV^{**}; u) \ \leq \ n_\cb(\OV;u).$$
\end{prf}

\medskip

\begin{rem}
\label{univ-cc}
(a) It is clear that if $\OV$ is a matricially normed space having a complete strict geometric unitary, then $\OV$ is an operator space (by Theorem \ref{ncb=1}). 
More generally, if $\OV$ has a complete geometric unitary, then $\OV$ is pseduo $L^\infty$ in the sense that there exists $\kappa \geq 1$ such that $\|u\oplus v\| \leq \kappa \max\{\|u\|, \|v\|\}$ for any $u\in M_m(V)$, $v\in M_n(V)$ and $m,n\in \mathbb{N}$ (by Corollary \ref{ncb+>ci}(a)). 
Consequently, if $p\in [1,\infty)$, then any $L^p$-matricially normed space will not have a complete geometric unitary. 

\smnoind
(b) Let $\OV$ be a matricially normed space and $u\in \sph(\OV)$. 
Suppose that $H$ is a Hilbert space and $\Theta:  \OV_u \rightarrow \CL(H)$ is a complete isometry such that $\Theta(Q_u(u)) = \id_H$ (by Lemma \ref{equiv-ncb>0}(b) and Theorem \ref{ncb=1}). 
Then $\Theta \circ Q_u:\OV \rightarrow \CL(H)$ satisfies certain universal property in the following sense: 
if $K$ is any Hilbert space and $\Psi: \OV \rightarrow \CL(K)$ is a complete contraction with $\Psi(u) = \id_K$, then there exists a complete contraction (not necessarily unique) $\Lambda: \CL(H) \rightarrow \CL(K)$ such that $\Psi = \Lambda \circ \Theta \circ Q_u$. 
Indeed, $\Psi = \Psi_u \circ Q_u$ (where $\Psi_u: \OV_u \rightarrow \CL(K)$ is as in Proposition \ref{map-ncb}(a)) and there exists, by the Averson extension theorem, a complete contraction $\Lambda: \CL(H) \rightarrow \CL(K)$ such that $\Psi_u = \Lambda \circ \Theta$. 
\end{rem}

\medskip

\subsection{Relationships with geometric unitaries}

\medskip

In this subsection, we will give some comparison between $n_\cb(\OV,u)$ and $n(V,u)$. 

\medskip

\begin{cor}
\label{ncb-n}
Let $\OV$ be a matricially normed space. 
If $n_\cb(\OV;u) > 0$ (respectively, $n_\cb(\OV;u) = 1$), then $n(V;u) > 0$ (respectively, $n(V;u)\geq \frac{1}{2}$). 
\end{cor}
\begin{prf}
Let $\Psi$ be the completely contractive complete topological injection (respectively, complete isometry) given by Corollary \ref{ncb+>ci}(a) (respectively, Theorem \ref{ncb=1}). 
Then $\Psi$ is a contractive topological injection (respectively, isometry) and Proposition \ref{map-n}(b) and \cite[Theorem 3]{CDM} shows that $n(V;u) \geq \frac{n(\CL(H);1)}{\|\Psi^{-1}\|} \geq \frac{1}{2\|\Psi^{-1}\|}$. 
\end{prf}

\medskip

Next, we compare $n$ and $n_\cb$ for the minimal quantization. 

\medskip

\begin{prop}
\label{n<ncbmin}
Let $X$ be a normed space and $u\in \sph(X)$. 

\smnoind
(a) $n(X;u) \leq n_\cb(\mathbf{\min X};u)$. 

\smnoind
(b) $n(X;u) > 0$ if and only if $n_\cb(\mathbf{\min X}; u)> 0$. 
\end{prop}
\begin{prf}
(a) Since $\CS_1(\mathbf{\min X}; u) = \CS(X;u)$, it suffices to prove that for any $k\in \mathbb{N}$ and $x = (x_{ij})\in \sph(M_k(\mathbf{\min X}))$, we have 
\begin{equation}
\label{n<sup}
n(X;u)\ \leq\ \sup \{\|f_k(x)\|: f\in \CS(X;u)\}.
\end{equation}
Suppose that $\Omega$ is a compact Hausdorff space such that $\mathbf{\min X} \subseteq C(\Omega)$ as operator subspace (hence, $M_k(\mathbf{\min X})\subseteq C(\Omega;M_k)$). 
There exist $\omega\in \Omega$ and $(c_i), (d_i)\in \sph(\mathbb{C}^n)$ with 
$$1\ =\ \|x\|\ =\ \|x(\omega)\|\ =\ \abs{\sum_{ij} \overline{c_i}x_{ij}(\omega)d_j}$$
which implies that $\|\sum_{ij} \overline{c_i}x_{ij}d_j\| \geq 1$. 
On the other hand, it is easy to see that $\|\sum_{ij} \overline{c_i}x_{ij}d_j\| \leq \|x\| = 1$.
Now, for any $f\in \CS(X;u)$, 
$$
\|f_k(x)\| 
\ \geq\ \abs{\left\langle \left(f(x_{ij})\right)
\left(\begin{array}{c}
d_1\\
\vdots\\
d_k
\end{array}\right), 
\left(\begin{array}{c}
c_1\\
\vdots\\
c_k
\end{array}\right)\right\rangle} 
\ = \ \abs{f\left(\sum_{ij} \overline{c_i}x_{ij}d_j\right)}.
$$
Hence, $\gamma^u(\sum_{ij} \overline{c_i}x_{ij}d_j) \leq \sup \{\|f_k(x)\|: f\in \CS(X;u)\}$ and Equality \eqref{n<sup} is verified. 

\smnoind
(b) This part follows from part (a) as well as Corollary \ref{ncb-n}. 
\end{prf}

\medskip

\begin{rem}
\label{rem-min}
(a) The arguement of the above actually shows that for any $(x_{ij}) \in M_k(\mathbf{\min X})$, there exists $(c_i), (d_i)\in \sph(\mathbb{C}^n)$ such that $\|(x_{ij})\| = \|\sum_{ij} \overline{c_i}x_{ij}d_j\|$. 
This could be a known fact although we do not find it in the literatures. 

\smnoind
(b) Let $X$ be a finite dimensional normed space and $\mathbf{X}$ be any quantization of $X$. 
Then the identity map is a completely contractive complete isomorphism from $\mathbf{X}$ to $\min \mathbf{X}$ (see e.g. \cite[2.2.4]{ER}). 
Consequently, by Propositions \ref{map-ncb}(b) and \ref{n<ncbmin}(b), if $u\in \sph(X)$, then $n_\cb(\mathbf{X},u) > 0$ if and only if $n(X,u) > 0$. 

\noindent
(c) If $\OV$ is an operator space and $u\in \sph(\OV)$ such that $n_\cb(\OV;u) > 0$, then $n_\cb(\mathbf{\min V};u) >0$ (by Corollary \ref{ncb-n} and Proposition \ref{n<ncbmin}).  
This fact is not easy to obtain directly from the definition. 
\end{rem}

\medskip

\subsection{An example}

\medskip

Let us start this subsection with the following result concerning $n_\cb$ for $C^*$-algebras. 
Note that (i)$\Rightarrow$(ii) follows from Theorem \ref{ncb=1} while (iii)$\Rightarrow$(i) follows from Corollary \ref{ncb-n} as well as \cite[4.1]{BK}. 

\medskip

\begin{cor}
For a $C^*$-algebra $A$ and $u\in \sph(A)$, the following statements are equivalent.

\smnoind
(i) $u$ is a unitary.

\smnoind
(ii) $n_\cb(A;u) = 1$. 

\smnoind
(iii) $n_\cb(A;u) > 0$. 
\end{cor}

\medskip

Note that from \cite[Theorem 3]{CDM}, if $A$ is a commutative unital $C^*$-algebra, then $n(A;1) = 1$ but if $A$ is noncommutative, then $n(A;1) = \frac{1}{2}$. 
By considering $n_\cb$, the above corresponding result looks cleaner but one losses the ability to detect whether $A$ is commutative. 

\medskip

On the other hand, for a general operator spaces $\OV$ and any $u,v\in \sph(\OV)$, it is possible that both $n_\cb(\OV,u)$ and $n_\cb(\OV,u)$ are non-zero but $n_\cb(\OV;u) \neq n_\cb(\OV;v)$. 
In order to give such an example, we need the following lemma which is probably known. 
However, since we cannot find it in the literatures, we give a proof here for completeness. 

\medskip

\begin{lem}
\label{n-sum}
Let $X$ and $Y$ be two normed spaces and $u\in \sph(X)$. 
Then $n(X\oplus^1 Y; (u,0)) = n(X;u)$. 
\end{lem}
\begin{prf}
Let $E = X\oplus^1 Y$. 
For any $(f,g) \in X^*\oplus^\infty Y^*$, we have $(f,g)\in \CS(E; (u,0))$ if and only if $\max\{\|f\|,\|g\|\} = 1 = f(u)$. 
Hence, $\CS(E; (u,0)) = \CS(X;u)\times {\rm Ball}(Y^*)$ (where ${\rm Ball}(Y^*)$ is the closed unit ball of $Y^*$). 
Thus, 
$$\gamma^{(u,0)}(x,y)\ \leq\ \gamma^u(x) +\|y\| \qquad ((x,y)\in E).$$
On the other hand, for any $\epsilon > 0$, there exists $f\in \CS(X;u)$ and $g\in {\rm Ball}(Y^*)$ such that $\gamma^u(x) < \abs{f(x)} + \epsilon$ and $g(y) = \|y\|$. 
If $f(x) = \abs{f(x)}e^{i\theta}$, then 
$\abs{f(x) + (e^{i\theta}g)(y)} = \abs{f(x)} + \|g\| \geq \gamma^u(x) + \|y\| - \epsilon$. 
This shows that $\gamma^{(u,0)}(x,y) = \gamma^u(x) + \|y\|$. 
Consequently, 
$$n(E;(u,0))\ =\ \inf \{ \gamma^{u}(x) +\|y\|: (x,y)\in E; \|x\| + \|y\| = 1\} \ \leq\ n(X;u)$$ 
(since $\|y\| = 0$ is possible). 
On the other hand, for any $n\in \mathbb{N}$, there exists $(x_n,y_n) \in E$ with 
$$\|x_n\| + \|y_n\|\ =\ 1 \quad {\rm and} \quad \gamma^u(x_n)+\|y_n\|\ <\ n(E;(u,0)) + \frac{1}{n}.$$ 
If there are infinitely many $n$ with $x_n = 0$, then there exists a subsequence such that $1 = \|y_{n_k}\| < n(E; (u,0)) + \frac{1}{n_k}$ which implies that $n(E;(u,0)) = 1$ and so $n(X;u) = 1$ as well. 
Otherwise, we can assume that all $x_n$ are non-zero and take $z_n = \frac{x_n}{1-\|y_n\|}\in \sph(X)$. 
Since 
$$(\gamma^u(x_n) + \|y_n\|)(1-\|y_n\|)\ =\ \gamma^u(x_n) + \|y_n\| (1 - \gamma^u(x_n) - \|y_n\|)\ \geq\ \gamma^u(x_n)$$ 
(note that $\gamma^u(x_n) \leq \|x_n\| = 1 - \|y_n\|$), we have $\gamma^u(z_n) \leq \gamma^u(x_n) + \|y_n\| < n(E;(u,0)) + \frac{1}{n}$. 
Therefore, $n(X;u) \leq n(E;(u,0))$ as required. 
\end{prf}

\medskip

\begin{eg}
\label{eg-nonuniq}
If $E$ is a finite dimensional Banach space with $n(E;u) = 1/e$ (see e.g. \cite[3.5]{DMPW}), then for any quantization $\mathbf{E}$ of $E$, we have $0 < n_\cb(\mathbf{E};u) < 1$ (by Remark \ref{rem-min}(b) and the fact that $1/e < 1/2$ together with Corollary \ref{ncb-n}). 
Moreover, if $F = \mathbb{C}\oplus^1 E$, then $n(F;(0,u)) = 1/e$ and $n(F; (1,0)) = 1$ (by Lemma \ref{n-sum}) and so, by the above, $0 < n_\cb(\mathbf{\min F}; (0,u)) < 1$ but $n_\cb(\mathbf{\min F}; (1,0)) = 1$ (by Proposition \ref{n<ncbmin}(a)). 
\end{eg}

\bigskip

\section{Applications}

\medskip

\subsection{An abstract definition for operator systems}

\medskip

If $\OV$ is an operator space, $u\in \sph(\OV)$ and $\Psi: \OV\rightarrow \CL(H)$ is a complete isometry such that $\Psi(u) = \id_H$, then it is known that $\Psi^{-1}(\Psi(\OV)\cap \Psi(\OV)^*)$ is an operator system which is independent of the choice of $(H, \Psi)$. 
The following result gives an alternative and intrinsic description for this operator system. 

\medskip

\begin{prop}
\label{os-sgu}
Suppose that $\OV$ is a matricially normed space and $u$ is a complete strict geometric unitary of $\OV$.  
For any $n\in \mathbb{N}$, we set 
$$\BK_u^n := \{x\in M_n(\OV): \varphi_n(x) \geq 0; \varphi \in \CS_k(\OV;u); k\in \mathbb{N}\}.$$
If $\Psi: \OV \rightarrow \CL(H)$ is any complete isometry with $\Psi(u) = \id_H$, then $\Psi_n(\BK_u^n) = \Psi_n(M_n(\OV))\cap \CL(H^{(n)})_+$.
Conseqently, $V_s = {\rm span}\ \BK_u^1$ is an operator system with matrix order structure given by $\BK_u^n\cap M_n(V_s)$. 
\end{prop}
\begin{prf}
If $A$ and $\Theta$ are as in the argument for Theorem \ref{ncb=1}, then clearly, $\Theta_n(\BK_u^n) = \Theta_n(M_n(\OV))\cap M_n(A)_+$. 
Regard $A$ as a unital $C^*$-subalgebra of some $\CL(K)$. 
If $\Psi: \OV \rightarrow \CL(H)$ is as in the statement, then there exist (by Remark \ref{univ-cc}(c)) complete contractions $\Gamma_1: \CL(H) \rightarrow \CL(K)$ and $\Gamma_2: \CL(K) \rightarrow \CL(H)$ such that $\Theta = \Gamma_1\circ \Psi$ and $\Psi = \Gamma_2\circ \Theta$. 
Hence, both $\Gamma_1$ and $\Gamma_2$ are completely positive, and so, $\Psi_n(\BK_u^n) = \Psi_n(M_n(\OV))\cap \CL(H^{(n)})_+$.
\end{prf}

\medskip

\begin{thm}
\label{ab-ch-os}
Let $\OV$ be a matrically normed space and $u\in \sph(V)$. 
Then there exists a matrix order structure on $\OV$ turning it into an operator system with order unit $u$ if and only if $n_\cb(\OV, u) = 1$ and $\BK_u^1$ spans $\OV$.
\end{thm}

\medskip

\begin{cor}
Let $\OV$ be matricially normed space and $u\in\sph(\OV)$. 

\smnoind
(a) $\OV_{u,s}$ is an operator system whose cone is $Q_u(\BK_u^1)$. 

\smnoind
(b) If $n_\cb(\OV;u) > 0$, there exists an equivalent matrix norm on $V$ under which $V_s$ is an operator system with order unit $u$. 

\smnoind
(c) If $\varphi:\OV\rightarrow \CL(H)$ is a complete contraction such that $\varphi(u) = \id_H$, then $\varphi$ factor through a ``completely positive map'' from $\OV_u$ to $\CL(H)$. 
\end{cor}

\medskip

One can also use the above ideal to give an abstract characterisation for ``non-unital operator systems''. 
However, since this is not directly related to complete geometric unitary, we will do this in the Appendix. 

\medskip

\subsection{Quotient by a complete $M$-ideal}

\medskip

In this subsection, we will show that the quotient of a unital operator space with a proper complete $M$-ideal is also unital. 
As an application, the quotient of any operator system by a proper complete $M$-ideal is again an operator system. 
Let us begin with the following lemma which may also be known. 

\medskip

\begin{lem}
\label{gu-sum}
Let $X$ and $Y$ be two normed spaces. 
If $(u,v)\in \sph(X\oplus^\infty Y)$ with $n(X\oplus^\infty Y; (u,v)) > 0$, then $\|u\| = 1 = \|v\|$. 
\end{lem}
\begin{prf}
Suppose that $\|u\| =1$ but $\|v\| < 1$. 
For any $(f,g)\in \CS(X\oplus^\infty Y; (u,v))$, we have $\|f\| + \|g\| = 1 = f(u) + g(v)$. 
Hence $g(v) = \|g\|$ and so, either $\|g\| = 0$ or $\|v\| \geq 1$ (which is absurd). 
Consequently, $\CS(X\oplus^\infty Y; (u,v)) = \CS(X; u) \times \{0\}$. 
Now, for any $y\in \sph(Y)$ and $(f,0)\in \CS(X\oplus^\infty Y; (u,v))$, we have $(f,0)(0,y) = 0$ and so, $n(X\oplus^\infty Y; (u,v)) = 0$ which is a contradiction.  
A similar contradiction occur if $\|v\| = 1$ but $\|u\| < 1$. 
\end{prf}

\medskip

If $A$ is a unital operator algebra with identity $1_A$ and $I$ is a proper closed ideal of $A$, then it is obvious that $1_A\notin I$ and the image of $1_A$ in $A/I$ is the identity. 
Interesting, this fact can be regarded as a geometric statement and can be generalized to the following result. 

\medskip

\begin{thm}
\label{quot-m-ideal}
(a) Let $\OV$ be an operator space and $\OW\subsetneq \OV$ be a complete $M$-ideal. 
If $v\in \sph(\OV)$ with $n_\cb(\OV;v) > 0$, then $\|Q(v)\| = 1$ and $n_\cb(\OV;v) \leq n_\cb(\OV/\OW; Q(v))$ (where $Q:\OV \rightarrow \OV/\OW$ is the canonical quotient map). 

\smnoind
(b) Let $\OV$ be an operator system and $\OW$ is a proper complete $M$-ideal of $\OV$. 
Then $\OV / \OW$ is also an operator system. 
\end{thm}
\begin{prf}
(a) By Corollary \ref{ncb+>ci}(b), one can assume that $\OW$ is a complete $M$-summand of $\OV$ (since $\OW^{\bot\bot}$ is a complete $M$-summand of $\OV^{**}$). 
Let $P: \OV\rightarrow \OV$ be the complete $M$-projection such that $P(\OV) = \OW$ and let $\OZ = (I-P)(\OV)$. 
Then $\OV \cong \OZ\oplus^\infty \OW$ as operator spaces and $Q(x) \mapsto (I-P)(x)$ is a complete isometry from $\OV/\OW$ to $\OZ$. 
Suppose that $v = (u,w)$ with $u\in \OZ$ and $w\in \OW$. 
Then by Corollary \ref{ncb-n} and Lemma \ref{gu-sum}, we see that $\|u\| = 1 = \|w\|$. 
Pick any $f\in \sph(\OZ^*)$ with $f(u) = 1$ and define $\Psi: \OZ \rightarrow \OV$ by $\Psi(x) = (x, f(x)w)$. 
It is not hard to check that $\Psi$ is a complete isometry (as $\|f_n(x)w\| \leq \|x\|$ for any $x\in M_n(\OZ)$) such that $\Psi(u) = (u,w)$. 
Now, $n_\cb(\OV; (u,w)) \leq n_\cb(\OZ;u)$ (by Proposition \ref{map-ncb}(b)). 

\smnoind
(b) Let $v$ be the order unit of $\OV$. 
Since the composition with the canonical quotient map $Q: \OV \rightarrow \OV / \OW$ is a map from $\CS(\OV / \OW; Q(v))$ to $\CS(\OV, v)$, we know that $Q(\BK_v^1)\subseteq \BK_{Q(v)}^1$ and so, $\BK_{Q(v)}^1$ spans $\OV / \OW$.
On the other hand, part (a) tells us that $n_\cb(\OV / \OW; Q(v)) =1$. 
Now, the result follows from Theorem \ref{ab-ch-os}. 
\end{prf}

\medskip

\begin{rem}
One can also show that the image of a (strict) geometric unitary of a Banach space in the quotient by a proper $M$-ideal is also a (strict) geometric unitary. 
\end{rem}

\bigskip

\appendix\section{Appendix: non-unital operator systems}

\bigskip

We begin this appendix with the following probably well-known result. 

\medskip

\begin{lem}
\label{pos-PTP}
Let $H$ be a Hilbert space and $T \in \CL(H^n)$. 
Then $T\geq 0$ if and only if $(P\otimes 1)T(P\otimes 1) \geq 0$ for any finite dimensional projection $P\in \CL(H)$. 
\end{lem}
\begin{prf}
We only need to show the sufficiency. 
Let $B$ be an orthonormal basis for $H$ and $H_0 := {\rm span}\ B$. 
Suppose that $\mathfrak{F}(B)$ is the collection of all finite subsets of $B$ and $P_F$ is the projection onto ${\rm span}\ F$ ($F\in \mathfrak{F}(B)$). 
Then $(P_F\otimes 1)T(P_F\otimes 1) \geq 0$ for any $F\in \mathfrak{F}(B)$ will imply that $\langle \eta, T \eta \rangle \geq 0$ for any $\eta\in H_0^n$, and so, $T \geq 0$. 
\end{prf}

\medskip

Let $\OV$ be a matrically normed space and $M_n(\OV)_+$ be a cone in $M_n(\OV)$ ($n\in \mathbb{N}$). 
For every $n\in \mathbb{N}$, we set 
$$\CS_n^+(\OV)\ :=\ \{\varphi\in \CB(\OV, M_n): \|\varphi\|_\cb \leq 1; \varphi_m (M_m(\OV)_+) \subseteq (M_{mn})_+; m\in \mathbb{N}\},$$ 
$$\gamma_k^+(x)\ :=\ \sup\left\{ \|\varphi_k(x)\|: \varphi\in \CS_n^+(\OV); n\in\mathbb{N}\right\} \qquad (k\in \mathbb{N}; x\in M_k (\OV)),$$
$$n_\cb^+(\OV)\ :=\ \inf\left\{\gamma_k^+(x): x\in \sph(M_k(\OV)); k\in \mathbb{N}\right\}$$
as well as
$$\BK_n\ :=\ \{v\in M_n(V): \varphi_n(v) \in (M_{nk})_+; k\in \mathbb{N}; \varphi\in \CS_k^+(\OV)\}.$$

\medskip

It is easy to check that $\CS_n^+(\OV)$ is compact under the point-norm topology and $M_n(\OV)_+ \subseteq \BK_n$ ($n\in \mathbb{N}$). 

\medskip

\begin{thm}
Suppose that $\OV$ and $M_n(\OV)_+$ be as in the above. 
Then there exist a Hilbert space $H$ and a complete isometry $\Phi: \OV \rightarrow \CL(H)$ with $\Phi(M_n(\OV)_+) = \Phi(M_n(\OV))\cap \CL(H^n)_+$ $(n\in \mathbb{N}$) if and only if 
$n_\cb^+(\OV) = 1$ and $M_n(\OV)_+ = \BK_n$ for any $n\in\mathbb{N}$.
\end{thm}
\begin{prf}
$\Rightarrow)$. 
$n_\cb^+(\OV) = 1$ because of Lemma \ref{supPTP} and $M_n(\OV)_+ = \BK_n$ follows from Lemma \ref{pos-PTP}. 

\smnoind
$\Leftarrow)$. 
Let $A := \bigoplus_{k=1}^\infty C(\CS_k^+(\OV), M_k)$ and define 
$\Theta: \OV \rightarrow A$ by $\Theta (v) = (\Theta^{(k)}(v))$ where $\Theta^{(k)}(v)(\varphi)= \varphi(v)$ ($\varphi\in V$). 
As in the argument for Theorem \ref{ncb=1}, $\Theta$ is a complete isometry (because $n_\cb^+(\OV) = 1$). 
Moreover, it is easy to see that $\Theta_n(\BK_n) = \Theta_n(M_n(\OV))\cap M_n(A)_+$ ($n\in \mathbb{N}$). 
\end{prf}

\medskip

This gives an abstract characterization of ``possibly non-self-adjoint and non-unital operator system'' as follows:
\begin{quotation}
Suppose that $\OV$ is a matrically normed space and $M_n(\OV)_+$ is a cone in $M_n(\OV)$ for any $n\in \mathbb{N}$. 
Then $(\OV, M_n(\OV)_+)$ is called an \emph{abstract (not necessarily unital) operator system} if $n_\cb^+(\OV) = 1$ and $M_n(\OV)_+ = \BK_n$ for any $n\in\mathbb{N}$. 
\end{quotation} 

\bigskip

\noindent{\bf Acknowledgement}

\smnoind
1. After this work has finished, D. Blecher and M. Neal recently gave another abstract characterization for unital operator spaces (as well as an operator space characterization for operator systems) in their paper \cite{BN}. 
More precisely, they characterize the ``identity'' $u$ of an operator space $X$ by looking at the norms of some matrices defined by $u$ together with an arbitrary elements in $M_n(X)$ ($n\in \mathbb{N}$). 

\smnoind
2. We would like to thank D. Blecher for some helpful comments on an earlier version of this work, especially for informing us the results in \cite{BRS}.

\medskip\noindent
Xu-Jian Huang, Chern Institute of Mathematics and LPMC, Nankai University, Tianjin 300071, China.

\smnoind
\emph{Email address:} huangxujian@mail.nankai.edu.cn

\medskip\noindent
Chi-Keung Ng, Chern Institute of Mathematics and LPMC, Nankai University, Tianjin 300071, China.

\smnoind
\emph{Email address:} ckng@nankai.edu.cn


\medskip



\bigskip

\begin{thebibliography}{99}

\bibitem{AW} C. Akemann and N. Weaver, Geometric characterizations of some classes of operators in $C\sp *$-algebras and von Neumann algebras, Proc. Amer. Math. Soc. 130 (2002), 3033-3037. 

\bibitem{ARS} O. Y. Aristov, V. Runde, Volker and N. Spronk, Operator biflatness of the Fourier algebra and approximate indicators for subgroups, J. Funct. Anal.  209 (2004), 367-387.

\bibitem{BL} D. Blecher and C. Le Merdy, \emph{Operator algebras and their modules---an operator space approach}, London Math. Soc. Mono. New Series 30, Oxford University Press (2004).

\bibitem{BN} D. Blecher and M. Neal, Metric characterizations of isometries and of unital operator spaces and systems, preprint (arXiv:0805.2166v1). 

\bibitem{BRS} D. Blecher, Z.-J. Ruan and A. Sinclair, A characterization of operator algebras, J. Funct. Anal. 89 (1990), 188-201. 

\bibitem{BK} H.F. Bohnenblust and S. Karlin, Geometrical properties of the unit sphere of Banach algebras, Ann. of Math. 62 (1955), 217-229.

\bibitem{BD} F.F. Bonsall and J. Duncan, \emph{Numerical ranges of operators on normed spaces and of elements of normed algebras}, London Math. Soc. Lect. Note Ser. 2, Cambridge University Press (1971).

\bibitem{CDM} M.J. Crabb, J. Duncan and C.M. McGregor, Characterizations of commutativity for $C\sp{*} $-algebras, Glasgow Math. J. 15 (1974), 172-175.

\bibitem{DMPW} J. Duncan, C.M. McGregor, J.D. Pryce and A.J. White, The numerical index of a normed space, J. London Math. Soc. 2 (1970), 481-488. 

\bibitem{ER} E.G. Effros and Z.-J. Ruan, \emph{Operator spaces}, London Math. Soc. Mono. New Series, 23, Oxford Univ. Press, (2000).

\bibitem{IS} M. Ilie and N. Spronk, Completely bounded homomorphisms of the Fourier algebras, J. Funct. Anal. 225 (2005), 480-499.

\bibitem{KR} J. Kraus and Z.-J. Ruan, Multipliers of Kac algebras, Internat. J. Math. 8 (1997), 213-248. 

\bibitem{LNR} A. Lambert, M. Neufang and V. Runde, Operator space structure and amenability for Fig\`{a}-Talamanca-Herz algebras, J. Funct. Anal. 211 (2004), 245-269.

\bibitem{LNW} C.W. Leung, C.K. Ng and N.C. Wong, Geometric unitaries in $JB$-algebras, preprint. 

\bibitem{Ng-cohom} C.K. Ng, Cohomology of Hopf $C\sp *$-algebras and Hopf von Neumann algebras, Proc. London Math. Soc. 83 (2001), 708-742.

\bibitem{Ng-ext} C.K. Ng, The Ext-functor for the category of completely bounded comodules, Internat. J. Math. 16 (2005), 307-332.

\bibitem{Ng-reg} C.K. Ng, Regular normed bimodules, J. Oper. Theory 56 (2006), 343-355. 

\bibitem{Pal} T.W. Palmer, \emph{Banach algebras and the general theory of $*$-algebras. Vol. 2. $*$-algebras}, Encyclopedia of Mathematics and its Applications 79, Cambridge University Press, Cambridge (2001). 

\bibitem{Ruan-char} Z.-J. Ruan, Subspaces of $C\sp *$-algebras, J. Funct. Anal. 76 (1988), 217-230. 

\bibitem{Ruan-op-amen} Z.-J. Ruan, The operator amenability of $A(G)$, Amer. J. Math. 117 (1995), 1449-1474.

\bibitem{Ruan-amen-Hopf} Z.-J. Ruan, Amenability of Hopf von Neumann algebras and Kac algebras, J. Funct. Anal. 139 (1996), 466-499.

\bibitem{Runde-os-aha} V. Runde, Applications of operator spaces to abstract harmonic analysis, Expo. Math. 22 (2004), 317-363.

\bibitem{Spr} N. Spronk, Operator weak amenability of the Fourier algebra, Proc. Amer. Math. Soc. 130 (2002), 3609-3617.
\end{thebibliography}
\end{document}